\documentclass[11pt]{article}
\usepackage{amssymb,latexsym,amsmath,amsbsy,amsthm,amsxtra,amsgen,dsfont,pdfsync,color}
\begin{document}
\title{A unifying approach to branching processes \\ in a varying environment}

\author{G\"otz Kersting\thanks{Institut f\"ur Mathematik, Goethe Universit\"at, Frankfurt am Main, Germany, kersting@math.uni-frankfurt.de, work partially supported by the DFG Priority Programme SPP 1590 ``Probabilistic Structures in Evolution''}}

\maketitle
\begin{abstract}
Branching processes $(Z_n)_{n \ge 0}$ in a varying environment  generalize the Galton-Watson process, in that they allow time-dependence of the offspring distribution. Our main results concern general criteria for a.s. extinction, square-integrability of the martingale $(Z_n/\mathbf E[Z_n])_{n \ge 0}$, properties of the martingale limit $W$ and  a Yaglom type result stating convergence to an exponential limit distribution of the suitably normalized population size $Z_n$, conditioned on the event $Z_n >0$. The theorems generalize/unify diverse results from the literature and lead to a classification of the processes.
\medskip


\end{abstract}

\section{Introduction and main results}

Branching processes $(Z_n)_{n \ge 0}$ in a varying environment  generalize the classical Galton-Watson processes, in that they allow time-dependence of the offspring distribution. This natural setting promises relevant applications (e.g. to random walks on trees as in \cite{Ru}) and received recently a renewal of interest, see e.g. \cite{Ba,BrHaut,GoKeMiPu,SaJa}. Former research on branching processes in a varying environment was temporarily affected by the appearence of certain exotic properties, and one could get the impression that it is difficult to grasp some kind of generic behaviour of these processes. Even so, steps in this direction were taken  by Peter Jagers \cite{Ja}, in particular, he aimed for a classification into  supercritical, critical and subcritical regimes in the spirit of ordinary Galton-Watson processes. In this paper we like to take up this line of research.
To this end we prove several theorems reaching from criteria for a.s. extinction  up to Yaglom type results. We  require only  mild regularity assumptions, in particular we  don't set any  restrictions to the  sequence of expectations $\mathbf E[Z_n]$, $n \ge 0$, thereby generalizing and unifying a number of individual results from the literature. 

In order to define a branching process in a varying environment (BPVE), let $Y_1, Y_2, \ldots$ denote a sequence of random variables with values in $\mathbb N_0$, and  $f_1,f_2, \ldots$ their distributions. Let $Y_{ni}$, $n,i\in \mathbb N$, be independent random variables such that $Y_{ni}$ and $Y_n$ coincide in distribution for all $n,i\ge 1$. Define the random variables $Z_n$, $n \ge 0$, with values in $\mathbb N_0$ recursively as
\[ Z_0 :=1 \ , \quad Z_{n} := \sum_{i=1}^{Z_{n-1}} Y_{ni} \ , \ n \ge 1\ . \]
Then the process $(Z_n)_{n \ge 0}$ is called a {\em branching process in the varying environment $v=(f_1,f_2, \ldots)$ with initial value $Z_0=1$.} These processes may be considered as a model for the development of the size of a population where individuals reproduce independently  with offspring distributions $f_n$ potentially changing among generations. Without further mention we always require that  $0 < \mathbf E[Y_n] < \infty$ for all $n \ge 1$.

There is one non-trivial statement on  BPVEs requiring no extra assumption. It says that $Z_n$ is a.s. convergent to a random variable $Z_\infty$ with values in $\mathbb N_0\cup \{\infty\}$. This result is due to Lindvall \cite{Li} and extends results of Church \cite{Chu} (for a comparatively short proof see  Theorem 1.4 in \cite{KeVa}). It also clarifies under which conditions $(Z_n)_{n \ge 0}$ may `fall asleep' at a positive state meaning that the event that $0< Z_\infty < \infty$ occurs with positive probability. Let us call such a branching process {\em asymptotically degenerate}. Thus for a BPVE it is no longer true that  the process  either gets extinct a.s. or else converges a.s. to infinity.  

As mentioned above a BPVE may exhibit  extraordinary properties, which don't show up  for ordinary Galton-Watson processes. Thus a BPVE may possess different growth rates, as detected by MacPhee and Schuh \cite{Schuh}. Here we establish a  framework which excludes such  exceptional  phenomena  and elucidates the generic behaviour. As we shall see, this is naturally done in an $\mathcal L^2$-setting. 

Our main assumption is a  uniformity requirement which reads as follows: There is a constant $c<\infty$ such that for all natural numbers $n \ge 1$ we have 
\begin{align} \label{A1}
\mathbf E[Y_n^2; Y_n \ge 2  ] \le c\, \mathbf E[Y_n; Y_n \ge 2] \cdot\mathbf E[Y_{n}  \mid Y_n \ge 1] \ .  \tag{A} 
\end{align}
This regularity assumption is considerably mild. As we shall explain in the next section, it is fulfilled for distributions $f_n$, $n \ge 1$, belonging to any common class of probability measures, like Poisson, binomial, hypergeometric, geometric, linear fractional, or negative binomial distributions, without  any restriction to the parameters.  It is  also satisfied    in the case that the random variables $Y_n$, $n \ge 1$, are a.s. uniformly bounded by a constant $c< \infty$.  To see this take into account that we have $\mathbf E[Y_{n}  \mid Y_n \ge 1] \ge 1$. 
Since a direct verification of \eqref{A1} may be tedious in examples, we shall present  in the next section a third moment condition which implies \eqref{A1} and which  can often  be easily  checked. 

Let us call a BPVE {\em regular}, if it fulfils condition \eqref{A1}.

\paragraph{Remark 1: A property of consistency.} Observe that together with a BPVE $(Z_n)_{n \ge 0}$ any subsequence $(Z_{n_i})_{i \ge 0}$ with $n_0:=0 < n_1 < n_2 <\ldots$ is a BPVE, too. We note that the condition \eqref{A1} is then transmitted, i.e. any subsequence of a regular BPVE is regular, too. The proof will be given after Lemma 6 below. \qed

\mbox{}\\
Before presenting our results let us agree on the following notational conventions: Let $\mathcal P$ be the set of all probability measures on $\mathbb N_0$. The weights of $f\in \mathcal P$ are named $f[k]$, $k \in \mathbb N_0$. We set
\[ f(s):= \sum_{k=0}^\infty s^k f[k] \ , \ 0 \le s \le 1\ . \]
Thus we denote the probability measure $f$ and its generating function by one and the same symbol. This facilitates presentation and will cause no confusion. Keep in mind that each operation applied to these measures has to be understood as an operation applied to their generating functions. Thus $f_1f_2$ stands not only  for the multiplication of  the generating functions $f_1,f_2$ but also for the convolution of the respective measures. Also $f_1 \circ f_2$ expresses the composition of  generating functions as well as the resulting probability measure. We shall consider the  mean and second factorial moment of a random variable $Y$ with distribution $f$,
\[ \mathbf E[Y]= f'(1)  \ , \ \mathbf E[Y(Y-1)] = f''(1) \ ,  \]

\noindent
and its normalized second factorial moment and normalized variance
\[ \nu:= \frac{\mathbf E[Y(Y-1)]}{\mathbf E[Y]^2}= \frac{f''(1)}{f'(1)^2} \ , \ \rho:= \frac{\mathbf{Var}[Y]}{\mathbf E[Y]^2}=\nu+ \frac 1{\mathbf E[Y]} -1\ .  \]

We shall discuss   branching processes in a varying environment  along the lines of  ordinary Galton-Watson processes. 
Let for $n \ge 1$
\[ q:= \mathbf P(Z_\infty =0) \ , \ \mu_n:= f_1'(1) \cdots f_n'(1) \ , \  \nu_n:= \frac{f_n''(1)}{f_n'(1)^2} \ , \  \rho_n:= \nu_n+ \frac 1{f_n'(1)}-1 \]
and also $\mu_0:=1$. Thus $q$ is the probability of extinction and $\mu_n= \mathbf E[Z_n]$, $n \ge 0$. Note that for the standardized factorial moments $ \nu_n$ we have  $ \nu_n < \infty$ under assumption \eqref{A1}. This implies $\mathbf E[Z_n^2] < \infty$  for all $n \ge 0$ (see Lemma 4 below).

Assumption \eqref{A1} is a mild requirement with substantial consequences, as seen from the following diverse necessary and sufficient criteria for a.s. extinction.

\paragraph{Theorem 1.} {\em Assume \eqref{A1}. Then the conditions
\begin{enumerate}
\item[\em (i)] $q=1$,
\item[\em (ii)] $\mathbf E[Z_n]^2 = o(\mathbf E[Z_n^2])^{\color{white}\big|}$ as $n \to \infty$,
\item[\em (iii)] $\displaystyle \sum_{k=1}^\infty \frac {\rho_k}{\mu_{k-1}}= \infty$, \vspace{-.35cm}
\item[\em (iv)] $\displaystyle \sum_{k=1}^\infty \frac {\nu_k}{\mu_{k-1}}= \infty$ \ or \ $ \mu_n \to 0$  
\end{enumerate}
are equivalent. Moreover, the conditions
\begin{enumerate}
\item[\em (v)] $q<1$,
\item[\em (vi)] $\mathbf E[Z_n^2] = O(\mathbf E[Z_n]^2)^{\color{white}\big|}$ as $n \to \infty$,
\item[\em (vii)] $\displaystyle \sum_{k=1}^\infty \frac {\rho_k}{\mu_{k-1}}< \infty$, \vspace{-.35cm}
\item[\em (viii)]  $\displaystyle \sum_{k=1}^\infty \frac {\nu_k}{\mu_{k-1}}< \infty$ \ and \ $\exists \, 0<r\le \infty: \mu_n \to r$  
\end{enumerate}
are equivalent.}

\newpage

\noindent
These conditions are useful in different ways.
Condition (iii)/(vii)  appears to be a particulary suitable  criterion for  a.s. extinction,  whereas the conditions (iv) and (viii)  will prove  helpful for the classification of BPVEs. Condition (vi) will allow us to determine the growth rate of $Z_{n}$, see Theorem 2. 
Observe that (ii) can be rewritten as $\mathbf E[Z_n]=o(\sqrt{\mathbf {Var}[Z_n]})$. Briefly speaking   this means that under \eqref{A1}  we have   a.s. extinction, if and only if the noise dominates the  average growth in the long run.

We point out that conditions (iii), (iv), (vii) and (viii)  access not only  the expectations $\mu_n$ but also the  second moments. This is a novel aspect in comparsion to ordinary Galton-Watson processes and also to Agresti's  classical  criterion  on BPVEs \cite[Theorem 2]{Agr}. Agresti's result provides a.s. extinction iff $\sum_{k\ge1} 1/\mu_{k-1}  =\infty$. He could do so by virtue of his  stronger assumptions, which exclude e.g. asymptotically degenerate processes. In our setting there is the possibility   that we have both $\sum_{k \ge 1} \rho_k/\mu_{k-1}=\infty$ and $\sum_{k \ge 1} 1/\mu_{k-1}<\infty$, and also  the other way round. This is shown by the following examples.

\paragraph{Example1.}  Let $Y_n$ take just the values $n+2$ and 0, with $\mathbf P(Y_n= n+2)= n^{-1}$. Then $\mathbf E[Y_n(Y_n-1)]\sim n$, $\mathbf E[Y_n]= 1+ 2/n$, $\mathbf E[Y_n-1 \mid Y_n\ge 1]\sim n$, thus \eqref{A1} is fulfilled. Also   $\mu_n \sim n^2/2$ and $\rho_n\sim n$, hence $\sum_{k \ge 1} 1/\mu_{k-1}<\infty$ and $\sum_{k \ge 1} \rho_k/\mu_{k-1}=\infty$. \qed

\paragraph{Example 2.} Let $Y_n$ take just the values 0,1 and 2, with $\mathbf P(Y_n=0)=\mathbf P(Y_n=2) = 1/(2n^2)$. Then $\mathbf E[Y_n(Y_n-1)] \sim n^{-2}$, $\mathbf E[Y_n]=1$ and $\mathbf E[Y_n-1\mid Y_n \ge 1]\sim 1/(2n^2)$, thus \eqref{A1} is fulfilled. Also $\mu_n=1$ and $\rho_n \sim n^{-2}$, hence $\sum_{k \ge 1} 1/\mu_{k-1}=\infty$ and $\sum_{k \ge 1} \rho_k/\mu_{k-1}<\infty$. \qed

\mbox{}\\
The last example exhibits an asymptotically degenerate branching process, as seen from the subsequent Corollary 1.

Next we  turn to the  normalized population sizes
\[ W_n := \frac {Z_n}{\mu_n} \ , \quad n \ge 0\ . \]
Clearly $(W_n)_{n \ge 0}$ constitutes a non-negative martingale, thus there exists an integrable random variable $W\ge 0$ such that we have
\[W_n \to W \text{ a.s.}\]
 as  $n \to \infty$. 
With \eqref{A1} the random variable $W$ exhibits the dichotomy known for  Galton-Watson processes. 

\paragraph{Theorem 2.} {\em For a regular BPVE we have:
\begin{enumerate}
\item[\em (i)]
If $q=1$, then $W=0$ a.s.
\item[\em (ii)]
If $q<1$, then $\mathbf E[W]=1$, $\mathbf E[W^2]<\infty$, and $\mathbf P(W=0)=q$.
\end{enumerate}}

\mbox{}\\
In particular, in case of $q<1$ the martingale $(W_n)_{n\ge 0}$ is convergent in $\mathcal L^2$ implying 
\begin{align}
\mathbf{Var}[W] =  \sum_{k=1}^\infty \frac {\rho_k}{\mu_{k-1}} \ .
\label{varW}
\end{align}
This formula goes back to Fearn  \cite{Fea}.  We point out that  Assumption \eqref{A1} excludes the possibility of $\mathbf P(W=0)>q$ and, in particular, of the possibility of different rates of growth as in the examples constructed by MacPhee and Schuh  \cite{Schuh} (see also \cite{Sou,Bi}).
By means of Theorem 2 (ii) we also gain further insight into  asymptotically degenerate processes. Under assumption \eqref{A1} they are just those processes which fulfil the properties $q<1$ and $0< \lim_{n\to \infty} \mu_n < \infty$. Also taking  Theorem 1 (v) and (viii) into account we obtain the following corollary.

\paragraph{Corollary 1.} {\em A regular BPVE is asymptotically degenerate, if and only if  both $\sum_{k=1}^\infty  \nu_k < \infty$ and the sequence $  (\mu_n)_{n \ge 0}$  has a positive, finite limit. Then $Z_\infty < \infty$ a.s.}

\mbox{}\\
Now we address  the  behaviour of the random variables $Z_n$ conditioned on the events that  \mbox{$Z_n>0$}. The next theorem shows that their values follow largely the corresponding conditional expectations $\mathbf E[Z_n \mid Z_n >0]$. For $n \ge 0$ let
\[a_n := 1+ \mu_n \sum_{k=1}^n \frac {\nu_k}{\mu_{k-1}}\ . \]

\paragraph{Theorem 3.} {\em For a regular BPVE,  the sequence of random variables $Z_n/a_n$ conditioned on $Z_n >0$, $n \ge 0$, is tight, i.e. for any $\varepsilon >0$ there is a $u<\infty$ such that for all $n \ge 0$
\begin{align} \mathbf P\Big( \frac{Z_n}{a_n} > u \mid Z_n >0\Big) \le \varepsilon \ , \label{bound1}
\end{align}
moreover, there exist numbers $\theta >0$ and $u>0$ such that for all $n\ge 0$
\begin{align} \mathbf P\Big( \frac{Z_n}{a_n} >u \mid Z_n >0\Big) \ge \theta \ . \label{bound2}
\end{align}
Also, we have
\begin{align} \gamma a_n \le \mathbf E[Z_n \mid Z_n >0] \le  a_n  
\label{EZn}
\end{align}
with some constant $\gamma >0$, so that we may replace $a_n$ by $\mathbf E[Z_n\mid Z_n>0]$  in \eqref{bound1} and \eqref{bound2}.}

\mbox{}\\
For $q<1$ we do not learn anything new from this theorem, here Theorem 2 (ii) gives much preciser information. Thus let us  focus on the case $q=1$, the situation of a.s. extinction. At first sight one might expect that the constant $\theta$  in \eqref{bound2} can be chosen arbitrarily close to 1, if only $u$ gets sufficiently small. This will  apply to many interesting cases, but it is not always true. The following example gives an illustration.

\paragraph{Example 3.} For $n \ge 1$ let
\[ f_{2n-1}[1]= 2^{-n} \ , \quad f_{2n-1}[0] =1-2^{-n} \quad \text{and} \quad f_{2n}[2^{n+1}-1]=f_{2n}[1] = \frac 12 . \]
It is easy to check that \eqref{A1} is valid (as well as the  conditions \eqref{B} and \eqref{C} below). 
We have $f_{2n-1}'(1)= 2^{-n}$ and $f_{2n}'(1)= 2^{n}$, hence
\[ \mu_{2n-1}= 2^{-n}  \quad \text{and} \quad \mu_{2n}= 1 \]
for all $n \ge 1$. In particular, we have $Z_{2n-1}\to 0$ in probability, which entails $q=1$. 
Also $\nu_{2n-1}=0$ and $\nu_{2n} \sim 2$ as $n \to \infty$ implying
\[    \sum_{k=1}^{2n} \frac{\nu_k}{\mu_{k-1}} \sim \sum_{k=1}^n 2^{k+1} \sim 2^{n+2} \quad \text{and}  \quad \sum_{k=1}^{2n-1} \frac{\nu_k}{\mu_{k-1}}= \sum_{k=1}^{2n-2} \frac{\nu_k}{\mu_{k-1}}\sim 2^{n+1}\ ,  \]
and
\[ a_{2n-1} \sim 3 \quad \text{and}Ê\quad a_{2n} \sim 2^{n+2} \ . \]
From Theorem 3 it  follows that there is a $z < \infty$ such that
\[ \mathbf P(Z_{2n-1} > z \mid Z_{2n-1} >0) \le \frac 12 \]
for all $n \ge 1$. Therefore
\begin{align*}
\mathbf P( Z_{2n} \le z \mid Z_{2n}>0) &= \mathbf P(Z_{2n}\le z \mid Z_{2n-1}>0) \notag \\&\ge \mathbf P(Z_{2n-1}\le z \mid Z_{2n-1} >0) f_{2n}[1]^z \notag \\&\ge 2^{-z-1} 
\end{align*}
for all $n \ge 1$, and for any $u>0$ 
\begin{align}\mathbf P( Z_{2n}/a_{2n} > u \mid Z_{2n}>0)  \le 1 - 2^{-z-1} 
\label{counterex} 
\end{align}
if $a_{2n} \ge z/u$. Since $a_{2n} \to \infty$, 
 the constant $\theta$ from \eqref{bound2} cannot take a value above $1 - 2^{-z-1}$ in this example.\qed

\mbox{}\\
This example suggests, that quite different scenarios may occur for BPVEs with $q=1$, and that their behaviour may abruptly change from one subsequence to the next. We point out that assumption \eqref{A1} does not put (e.g. for Poisson distributions) any restrictions onto the expectation $\mu_n$, $n \ge 1$, allowing a variety of examples. Of special interest is the case that the numbers $a_n$ are uniformly bounded. Here Theorem 3 reads as follows.

\paragraph{Corollary 2.} {\em Under assumption \eqref{A1} the conditions
\begin{enumerate}
\item[\em (i)] the sequence of random variables $Z_n$ conditioned on the events that $Z_n >0$, $n \ge 0$, is tight, 
\item[\em (ii)] $\sup_{n\ge 0} \mathbf E[Z_n \mid Z_n >0] < \infty$,
\item[\em (iii)] $\displaystyle \sum_{k=1}^n \frac{\nu_k}{\mu_{k-1}} = O \Big( \frac 1{\mu_n}\Big)$ as $n \to \infty$,
\end{enumerate}
are  equivalent.}

\mbox{}\\
For   an ordinary Galton-Watson process   these three conditions apply  just  in the subcritical regime, then  the conditioned random variables $Z_n$  have even a limiting distribution. It is easy to see that such a  feature will not hold in  general for a BPVE. Indeed: there are two offspring distributions $\hat f$ and $\tilde f$  such that the  limiting distributions $\hat g$ and $\tilde g$ for the corresponding conditional Galton-Watson processes differ from each other. Choose an increasing sequence $0=n_0< n_1 <n_2 < \cdots$ of natural numbers and  consider the BPVE $(Z_n)_{n \ge 0}$ in the varying environment $v=(f_1,f_2, \ldots)$, where $f_n= \hat f$ for $n_{2k} < n \le n_{2k+1}$, $k \in \mathbb N_0$, and $f_n = \tilde f$ else. Then it is obvious that  $Z_{n_{2k+1}}$ given the event $Z_{n_{2k+1}}>0$ converges in distribution to $\hat g$ and $Z_{n_{2k}}$ given the event $Z_{n_{2k}}>0$ converges in distribution to $\tilde g$,  provided that the sequence $(n_k)_{k \ge 0}$ is increasing sufficiently fast.

Thus it may come as a surprise that in the opposite situation of $1/\mu_{n}=o\big(\sum_{k=1}^{n} \nu_k/\mu_{k-1}\big)$ we encounter a distinctive behaviour of the conditional limit distributions of $Z_n$, which is in accordance with Yaglom's theorem  for ordinary Galton-Watson processes. For technical reasons we have to somewhat strengthen assumption \eqref{A1}. We require that for every $\varepsilon >0$ there is a constant $c_\varepsilon < \infty$ such that for all natural numbers $n \ge 1$
\begin{align} \label{B} 
\mathbf E\big [ Y_n^2 ; Y_n > c_\varepsilon(1+ \mathbf E[Y_n])\big] \le \varepsilon \mathbf E \big[ Y_n^2; Y_n \ge 2\big] 
\tag{B} 
\end{align}
This condition is again widely satisfied, as we shall explain in the next section.  It implies assumption \eqref{A1}. Namely, for $\varepsilon=1/2$ we have
\begin{align}\mathbf E[Y_n^2; Y_n \ge 2] &\le 2 \mathbf E[ Y_n^2 ; 2 \le Y_n \le c_{1/2}(1+ \mathbf E[Y_n])] \le 2c_{1/2} (1+ \mathbf E[Y_n]) \mathbf E[Y_n; Y_n \ge 2] \ .
\label{BimpliesA}
\end{align}
Since $1+ \mathbf E[Y_n] \le 2\mathbf E[Y_n \mid Y_n \ge 1]$, we obtain \eqref{A1} with $c=4c_{1/2}$. 

\paragraph{Theorem 4.} {\em Let \eqref{B} be satisfied and let $q=1$. Then the following conditions are equivalent:
\begin{enumerate}
\item[\em (i)] There is a sequence $b_n$, $n \ge 0$, of positive numbers such that   $Z_n/b_n$ conditioned on the event $Z_n>0$ converges in distribution to a standard exponential distribution as $n \to \infty$,
\item[\em (ii)] $\mathbf E[Z_n \mid Z_n >0] \to \infty$ as $n \to \infty$,
\item[\em (iii)] $\displaystyle  \frac1{\mu_n} = o \Big( \sum_{k=1}^n \frac {\nu_k}{\mu_{k-1}} \Big)$ as $n \to \infty$.
\end{enumerate}
Under these conditions we may set $b_n:= \mathbf E[Z_n \mid Z_n >0]$, and we have
\[ \mathbf E[Z_n \mid Z_n >0] \sim   \frac{\mu_n }2\sum_{k=1}^n \frac {\nu_k}{\mu_{k-1}} \ ,\]
or equivalently
\[ \mathbf P(Z_n>0) \sim 2\Big(\sum_{k=1}^n \frac {\nu_k}{\mu_{k-1}}\Big)^{-1}  \]
as $n \to \infty$.}

\mbox{}\\
This theorem covers the classical results of Kolmogorov  and Yaglom for critical Galton-Watson processes in the finite variance case (without further moment restrictions), since then \eqref{B} is trivially satisfied. 

\mbox{}\\
Our results show the way how to implement a classification of regular BPVEs, which connects to the  notions used for classical Galton-Watson processes. If $q<1$, then in view of Theorem~2 and  Corollary 1 we distinguish two regimes. There is the {\em supercritical regime} in the case of $\mathbf E[Z_n]\to \infty$, and the {\em asymptotically degenerate regime} otherwise. If on the other hand we have $q=1$, then Theorem 4 suggests to characterize the {\em critical regime}  by the condition $\mathbf E[Z_n \mid Z_n>0]\to \infty$ (and not by just   some condition on the limiting behaviour of $\mu_n$, as one might do in a first attempt), and  to allocate  the other BPVEs to the {\em subcritical regime}. In this way we  differentiate  the clear-cut  limiting property of critical BPVEs from the indeterminacy of the remaining processes. In this classification a subcritical BPVE $(Z_n)_{n\ge 0}$  exhibits subcritical behaviour in the sense that according to Theorem 3  the random variables $Z_n$ conditioned on $Z_n>0$ are tight at least along some subsequence, in which the $a_n$ stay bounded.  The $Z_n$  may   diverge with positive probability along some other  subsequence, yet this does  in general not imply critical behaviour in the sense that along that subsequence the random variables $Z_{n}$, conditioned on  $Z_{n}>0$ and suitably scaled, have  asymptotically an exponential distribution. For a counter-example we refer to the construction leading to   formula \eqref{counterex}.

By means of  Theorem 1 and Theorem 3 we may streamline the determining conditions of the four regimes, as summerized in the subsequent overview.

\paragraph{Proposition 1.} {\em A regular BPVE is
\begin{align*}
\text{{ supercritical}, iff\ } & \quad \lim_{n\to \infty} \mu_n=\infty\  \text{and} \ \sum_{k=1}^\infty \frac{\nu_k}{\mu_{k-1}} < \infty\ ,\\
\text{{ asymptotically degenerate}, iff\ } & \quad0< \lim_{n\to \infty} \mu_n<\infty \ \text{and}\ \sum_{k=1}^\infty \frac{\nu_k}{\mu_{k-1}} < \infty\ ,\\
\text{{ critical}, iff\ } & \quad   \lim_{n \to \infty} \mu_n \sum_{k=1}^n \frac {\nu_k}{\mu_{k-1}} =\infty \ \text{and}\ \sum_{k=1}^\infty \frac{\nu_k}{\mu_{k-1}} = \infty \ ,\\
\text{{ subcritical}, iff\ }&\quad\liminf_{n \to \infty}\mu_n=0 \ \text{and}\ \liminf_{n\to \infty}\mu_n\sum_{k=1}^n \frac {\nu_k}{\mu_{k-1}} <\infty \ .
\end{align*}}

\noindent
Note that convergence of the means $\mu_n$ is not enforced in the critical case, they may diverge, converge to zero or even oszillate in between.

\paragraph{Example 5.} 
 In the case $0 < \inf_n \nu_n \le \sup_n \nu_n < \infty$ (as e.g. for Poisson variables) the classification simplifies. Here we are in the supercritical regime, iff $\sum_{k \ge 0} 1/\mu_k < \infty$ (enforcing $\mu_n \to \infty$). Asymptotically degenerate behaviour is excluded, and there is plenty of room for critical processes, i.e. for  processes which conform to the conditions  $\sum_{k \ge 0} 1/\mu_k = \infty$ and $1/\mu_n =o\big(\sum_{k=0}^{n-1} 1/\mu_k\big)$. The second requirement is e.g. fulfilled, if we have $\mu_n/\mu_{n-1} \to 1$ as $n \to \infty$. This latter condition covers a variety of scenarios for $\mu_n$ below exponential growth and above exponential decay. \qed

\paragraph{Example 6.} In the {\em binary case} $\mathbf P(Y_n=2) = p_n$, $\mathbf P(Y_n=0)=1-p_n$ we get $f_n'(1)= f_n''(1)=2p_n$.  Therefore $\nu_k/\mu_{k-1}=1/\mu_k$, so that  the situation conforms to  the previous example. \qed

\paragraph{Example 7.} In the {\em symmetric case} $\mathbf P(Y_n=0)=\mathbf P(Y_n=2)= p_n/2$ and $\mathbf P(Y_n=1)=1-p_n$ we have $\mu_n=1$ and $\nu_n=p_n$. Here we find critical or asymptotically degenerate behaviour, according to whether $\sum_{k=1}^\infty p_n$ is divergent or convergent. \qed

\paragraph{Example 8.} If the $Y_n$ take only the values 0 and 1, then all $\nu_n$ vanish. Now the BPVE is subcritical or asymptotically degenerate, according to whether $\mu_n$ converges to zero or to a positive value.\qed

\noindent
Our proofs rely largely on analytic considerations. The task is to get a grip on the probability measures $f_1 \circ \cdots \circ f_n$, which are the distributions of the random variables $Z_n$. In order to handle such iterated compositions of generating functions we resort to a device which has been applied from the beginning in the theory of branching processes. For a probability distribution $f$ on $\mathbb N_0$ with positive, finite mean $m$ we define a function $\varphi:[0,1)\to \mathbb R$ by the equation
\[  \frac 1{1-f(s)} = \frac 1{m(1-s)} + \varphi(s) \ , \quad 0 \le s < 1 \ . \]
In this way the mean and the `shape' of $f$ are separated to a certain extent. Indeed,  Lemma 1 below shows that $\varphi$ takes values which are of the size of the {\em standardized} second factorial moment $\nu$. Therefore we briefly name $\varphi$ the {\em shape function} of $f$.
As we shall see these functions are useful to dissolve the generating function $f_1 \circ \cdots \circ f_n$ into a sum (see Lemma 5 below). Here our contribution consists in obtaining sharp upper and lower bounds for the function $\varphi$ and its derivative. The interaction of these bounds then allows for   precise estimates  e.g. of the survival probabilities $\mathbf P(Z_n >0)$. The role of assumption \eqref{A1} in this interplay is to keep both bounds together uniformly in $n$.

Concluding this introduction let us comment on the literature. Agresti in his paper \cite{Agr}  on a.s. extinction already derived the sharp upper bound for $\varphi$ which we  give below in formula \eqref{phi}. We note that this bound is related to the well-known Paley-Zygmund inequality (compare the proof of Lemma 7). Agresti also obtained a lower bound for the survival probabilities, which, however, in general is away from our sharp bound.   Lyons \cite{Ru}  obtained the equivalence of  the conditions (v), (vi), (vii) and (somewhat disguised) (viii) from Theorem 1 under the assumption that the random variables $Y_n$ are a.s. bounded by a constant, with methods  completely different from ours. He also proved Theorem 2, again under the assumption that the offspring numbers are a.s. uniformly bounded  by a constant. D'Souza and Biggins \cite{Bi} derived Theorem 2 under a different set of assumptions. They require that there are numbers $a>0, b>1$ such that $\mu_{m+n}/\mu_m \ge ab^n$ for all $m,n \ge 1$ (called the uniform supercritical case). They do not need finite second moments but assume instead that the random variables $Y_n$ are uniformly dominated by a random variable~$Y$ with  $\mathbf E[Y\log^+ Y]< \infty$. Goettge \cite{Goe} obtains $\mathbf E [W]=1$ under the  condition $\mu_n \ge an^b$ with $a>0,b>1$ (together with a uniform domination assumption), but  doesn't consider the validity of the equation $\mathbf P(W=0)=q$.   In order to prove the conditional limit law   from Theorem~4 Jagers \cite{Ja} draws attention to uniform estimates due to Sevast'yanov \cite{Sev} (see also Lemma 3 in~\cite{Fah}). This approach demands amongst others the strong  assumption that the sequence $\mathbf E[Z_n]$, $n \ge0$, is bounded from above and  away from zero.  Independently and in parallel to our work Bhattacharya and Perlman \cite{Bha}  have presented a considerable generalization of Jager's result,  on a different route and under assumptions which are stronger  than ours. For recent results on a.s. extinction and asymptotic exponentiality  of multitype BPVEs we refer to~\cite{Do}.

The paper is organized as follows. In Section 2 we discuss the assumptions and several examples. In Section 3 we analyze the shape function $\varphi$.  Section 4 contains the proofs of our theorems.  

\newpage

\section{Examples}

The following  example illustrates the difference in range of the conditions \eqref{A1} and \eqref{B}.

\paragraph{Example 9.} Let $Y$ have a {\em linear fractional distribution} meaning that
\[ \mathbf P(Y=y \mid Y\ge 1)= (1-p)^{y-1}p  \ , \ y \ge 1 \]
with some $0<p<1$ and some probability $\mathbf P(Y\ge 1)$. Then from properties of geometric distributions we have
\begin{align*}  
\mathbf E[Y\mid Y \ge 1]= \frac 1p \ , \ \mathbf E[Y-1\mid Y \ge 1] =\frac {(1-p)}p\ , \ \mathbf E[Y(Y-1)\mid Y\ge 1]= \frac{2(1-p)}{p^2} \ ,  \end{align*}
and it follows 
\begin{align*} 
\mathbf E[ Y^2; Y \ge 2 ] &\le 2 \mathbf E[Y(Y-1)]=\frac{4(1-p)}{p^2} \mathbf P(Y \ge 1) \\&= 4 \, \mathbf E[Y-1;Y \ge 1] \cdot \mathbf E[Y\mid Y \ge 1]  \le 4\, \mathbf E[(Y;Y \ge 2] \cdot \mathbf E[Y\mid Y \ge 1] \ .
\end{align*}
Thus for {\em any} sequence $Y_n$ of linear fractional random variables  assumption \eqref{A1} is fulfilled with $c=4$, whatever their parameters $p_n$ and $\mathbf P(Y_n \ge 1)$ are. 

However, for condition \eqref{B} the corresponding statement fails. To see this we resort for linear fractional distributions to the formula
\[\frac{2(1-p_n)}{p_n^2} \mathbf P(Y_n\ge 1)=\mathbf E[Y_n(Y_n-1)] \le \mathbf E[Y_n^2; Y_n \ge 2] \ .\]
If we assume \eqref{B}, then   also the inequality \eqref{BimpliesA} is valid yielding
\[\frac{2(1-p_n)}{p_n^2} \mathbf P(Y_n\ge 1) \le 4 c_{1/2}  \mathbf E[ Y_n-1; Y_n\ge 1] \cdot (1+ \mathbf E[Y_n]) \ . \]
For linear fractional distributions this estimate may be rewritten as
\[\frac{2(1-p_n)}{p_n^2} \mathbf P(Y_n\ge 1) \le 4 c_{1/2} \frac {(1-p_n)}{p_n}\mathbf P(Y_n \ge 1)\Big(1+ \frac 1p_n \mathbf P(Y_n \ge 1)\Big)\ , \]
which simplifies to
\[ \frac 1{2c_{1/2}} \le p_n + \mathbf P(Y_n \ge 1) \ . \]
Thus  condition \eqref{B} implies $\inf_n (p_n + \mathbf P(Y_n \ge 1)) >0$, and a sequence of linear fractional random variables satisfying $p_n + \mathbf P(Y_n \ge 1)\le1/n$ does not meet \eqref{B}. 

Incidentally, Theorem 4 still holds true for  linear fractional $Y_n$, $n \ge 1$,  regardless of the validity of \eqref{B}. Then, as is well known, also $Z_n$ is linear fractional for any $n \ge 1$, and consequently  the sequence $Z_n/ \mathbf E[Z_n \mid Z_n\ge 1]$ given the events that $Z_n \ge 1$ converges in distribution to a standard exponential distribution provided that we have $\mathbf E[Z_n \mid Z_n\ge 1] \to \infty$. 
\qed

\mbox{}\\
In other examples a direct verification of assumptions \eqref{A1} or \eqref{B} can be cumbersome. Therefore we introduce another assumption, which often is easier to handle. It reads: There is a constant $\bar c < \infty$ such that for all natural numbers $n \ge 1$
\begin{align} \mathbf E[Y_n(Y_n-1)(Y_n-2)] \le \bar c\, \mathbf E[Y_n(Y_n-1)] \cdot(1+ \mathbf E[Y_n]) \tag{C}
\label{C}
\end{align}  
Condition \eqref{C} implies \eqref{A1} and \eqref{B}, as seen from the following proposition.

\paragraph{Proposition 2.}{\em If condition \eqref{C} is fulfilled, then \eqref{B} holds with $c_\varepsilon :=\max(3,5\bar c/\varepsilon)$ and \eqref{A1} holds with $c:=\max (12,40 \bar c)$.}

\begin{proof}
From $c_\varepsilon \ge 3$ and \eqref{C} we obtain
\begin{align*}
\mathbf E[ Y_n^2; Y_n > c_\varepsilon(1+ \mathbf E[Y_n])] &\le 5\, \mathbf E[(Y_n-1)(Y_n-2); Y_n > c_\varepsilon(1+ \mathbf E[Y_n])] \\
&\le 5\, \frac {\mathbf E[Y_n(Y_n-1)(Y_n-2)]}{c_\varepsilon(1+ \mathbf E[Y_n])} \\
&\le \frac {5\bar c}{c_\varepsilon} \, \mathbf E[Y_n(Y_n-1)] \ .
\end{align*}
It follows
\[ \mathbf E[ Y_n^2; Y_n > c_\varepsilon(1+ \mathbf E[Y_n])] \le \varepsilon \, \mathbf E[Y_n^2; Y_n \ge 2]\ ,\]
which is our first claim. The second one follows by means of  \eqref{BimpliesA}.
\end{proof}

\mbox{}\\
Condition \eqref{C} can be easily handled by means of generating functions and its derivatives. Here are some examples.

\paragraph{Example 10.}  If the $Y_n$ are a.s. uniformly bounded by a constant $c$, then \eqref{C} is satisfied with $\bar c=c$. \qed

\paragraph{Example 11.} Let $Y$ be {\em Poisson} with parameter $\lambda >0$. Then 
\[ \mathbf E[Y(Y-1)(Y-2)] = \lambda^3 \le \lambda^2(\lambda +1)= \mathbf E[Y(Y-1)](1+\mathbf E[Y]) \ . \]
Here \eqref{C} is fulfilled with $\bar c=1$. \qed

\paragraph{Example 12.} For {\em binomial} $Y$  with parameters $m\ge 1$ and $0 < p < 1$ the situation is analog, here
\begin{align} \mathbf E[Y(Y-1)(Y-2)]= m(m-1)(m-2)p^3  \le  m(m-1)p^2mp \le\mathbf E[Y(Y-1)](1+ \mathbf E[Y]) \ .  \tag*{$\Box$}
\end{align}

\paragraph{Example 13.} For  a {\em hypergeometric distribution} with parameter $(N,K,m)$    we have for $N \ge 3$
\begin{align*}\mathbf E[Y(Y-1)(Y-2)]&=\frac{m(m-1)(m-2)K(K-1)(K-2)}{N(N-1)(N-2)} \\
&\le 3\frac{m(m-1) K(K-1) }{N(N-1) }\frac{mK}N \le 3  \mathbf E[Y(Y-1)](1+\mathbf E[Y])\ ,
\end{align*} 
and \eqref{C} is satisfied with $\bar c=3$. The case $N \le 2$ can  immediately be included.

(v) For {\em negative binomial distributions}  the generating function is given by
\[ f(s)= \Big( \frac p{1-s(1-p)}\Big)^\alpha \]
with $0<p<1$ and a positive integer $\alpha$. Now
\begin{align*}\mathbf E[Y]= \alpha \frac {1-p}p \ , \ \mathbf E[Y(Y-1)] = \alpha(\alpha+1) \frac {(1-p)^2}{p^2} \ , \\ 
\mathbf E[Y(Y-1)(Y-2)] = \alpha(\alpha +1)(\alpha+2) \frac {(1-p)^3}{p^3}\ . 
\end{align*}
Thus
\begin{align*}
\mathbf E[Y(Y-1)(Y-2)] \le 3 \mathbf E[Y(Y-1)](1+ \mathbf E[Y]) \ .
\end{align*}
Again \eqref{C} is fulfilled with $\bar c=3$. \qed

\section{Bounds for the shape function}

For $f \in \mathcal P$ with mean $0<m=f'(1)<\infty$ define the {\em shape function} $\varphi=\varphi_f: [0,1) \to \mathbb R$ by the equation
\[ \frac 1{1-f(s)} = \frac{1}{m(1-s)} + \varphi(s) \ , \ 0\le s < 1 \ . \]
Due to convexity of $f(s)$ we have $\varphi(s) \ge 0$ for all $0\le s < 1$.
By means of a Taylor expansion of $f$ around 1 one obtains
$\lim_{s \uparrow 1} \varphi(s) =   f''(1)/(2f'(1)^2) $,
thus we extend $\varphi$ by setting
\begin{align} \varphi(1):=\frac \nu 2 \quad \text{with}\quad \nu:=   \frac{f''(1)}{f'(1)^2}  \ . \label{phi1}
\end{align}
In this section we prove the following sharp bounds. 

\paragraph{Lemma 1.} {\em Assume $f''(1)<\infty$. Then for $0\le s \le 1$}
\begin{align}\label{phi}  \frac 12 \varphi(0) \le \varphi(s) \le 2 \varphi(1) \ . \end{align}

\mbox{}\\
Note that $\varphi$ is identical zero  if $f[z]=0$ for all $z \ge 2$. Else   $\varphi(0)>0$, and  the lower bound of $\varphi$ becomes strictly positive. Choosing $s=1$ and $s=0$ in \eqref{phi} we obtain $\varphi(0)/2\le \varphi(1)$ and $\varphi(0)\le 2\varphi(1)$. Note  that for $f=\delta_k$ (Dirac-measure at point $k$) and $k \ge 2$ we have $\varphi(1)=\varphi(0)/2$ implying that the constants 1/2 and 2 in \eqref{phi} cannot be improved. The upper bound was derived in \cite{Gei} using a different method of proof.

The next lemma is based on a close investigation of the derivative of $\varphi(s)$. 

\paragraph{Lemma 2.} {\em Let  $Y$ be a random variable with distribution $f$ and assume $ f''(1) < \infty$. Then for $0 \le s \le 1$ and  natural numbers $a\ge 1$}
\[ \sup_{s \le t \le 1}| \varphi(1)-\varphi(t)| \le 2 m\nu^2 (1-s)+2a\nu (1-s)+ \frac 2{m^2}\mathbf E[Y^2; Y > a] \ . \]

\mbox{}\\
Uniform estimates of $\varphi(1)-\varphi(s)$ based on third moments have already been obtained by Sevast'yanov \cite{Sev} and others (see Lemma 3 in \cite{Fah}). Our lemma implies and generalizes these estimates. For the proof of these lemmas we use the following  result.

\paragraph{Lemma 3.} 
{\em Let $g_1,g_2$ be elements of $\mathcal P $ with the same support and satisfying the following property: For any $y\in \mathbb N_0$ with $g_1[y]>0$ we have
\[  \frac {g_1[z]}{g_1[y]} \le \frac {g_2[z]}{g_2[y]}  \text{ for all } z>y \ . \]
Also let $\alpha: \mathbb N_0 \to \mathbb R$ be a non-decreasing function. 
Then}
\[ \sum_{y=0}^\infty \alpha(y) g_1[y] \le \sum_{y=0}^\infty \alpha(y) g_2[y] \ . \]

\begin{proof} The lemma's assumption is called  the `monotone likelihood ratio property', which  is known to imply our claim. For convenience, we give a short proof:
By assumption there is a  non-decreasing function $h(y)$, $y \in \mathbb N_0$, such that $h(y)= g_2(y)/g_1(y)$ for all elements $y$ of the support of $g_1$. Then for any real number $c$
\begin{align*}
 \sum_{y=0}^\infty \alpha(y) g_2[y]- \sum_{y=0}^\infty \alpha(y) g_1[y]  = \sum_{y=0}^\infty (\alpha(y)-c)(g_2[y]-g_1[y])  
 = \sum_{y=0}^\infty (\alpha(y)-c)(h(y)-1) g_1[y] \ .
\end{align*}
For $c:= \min\{ \alpha(y): h(y) \ge 1\}$ we have $\alpha(0) \le c < \infty$. For this choice of $c$, since $h$ and $\alpha$ are non-decreasing, every summand of the right-hand sum is non-negative. Thus the whole sum is non-negative, too, and our assertion follows.
\end{proof}

\begin{proof}[Proof of Lemma 1]
(i) First we examine a special case of Lemma 3. Consider for $0 < s \le 1$ and $r \in \mathbb N_0$ the probability measures
\begin{align*}   g_s[y]=  \frac {s^{r-y}}{1+ s+ \cdots+ s^{r}} \ , \quad  0\le y \le r \ .\end{align*}
Then for $0 < s \le t\le 1$, $0 \le y<z\le r$ we have $g_s[z ]/g_s[y]=s^{y-z}\ge t^{y-z}=   g_t[z]/g_t[y]$. Hence we obtain that
\begin{align*}  \sum_{y=0}^r yg_s[y]= \frac{s^{r-1}+2 s^{r-2} + \cdots + r}{1+ s+ \cdots+ s^{r}}
\end{align*}
is a decreasing function in $s$. Also $\sum_{y=0}^r yg_0[y] = r$ and $\sum_{y=0}^r yg_1[y]= r/2$, and it follows for $0\le s \le 1$
\begin{align}\label{upperlower} \frac r2  \le \frac{r+(r-1)s+ \cdots + s^{r-1}}{1+ s + \cdots + s^{r}} \le r\ .
\end{align}

(ii) Next we derive a second representation for $\varphi$. We have
\[ 1- f(s)= \sum_{z=1}^\infty f[z] (1-s^z) = (1-s) \sum_{z=1}^\infty f[z] \sum_{k=0}^{z-1} s^k\ ,\]
and
\begin{align*} f'(1)(1-s)- (1-f(s))&= (1-s) \sum_{z=1}^\infty f [z] \sum_{k=0}^{z-1} (1-s^k) \\
& = (1-s)^2 \sum_{z=1}^\infty f [z] \sum_{k=1}^{z-1} \sum_{j=0}^{k-1} s^j \\
&= (1-s)^2 \sum_{z=1}^\infty f [z] ((z-1) + (z-2)s  + \cdots + s^{z-2})\ .
\end{align*}
Therefore
\begin{align*} \varphi(s) &= \frac{m(1-s)- (1-f(s))}{m(1-s)(1-f(s))} \\
&= \frac{\sum_{y=1}^\infty f [y] ((y-1) + (y-2)s  + \cdots + s^{y-2})}{m\cdot \sum_{z=1}^\infty f [z] (1+ s  + \cdots + s^{z-1})} \ . 
\end{align*}
From \eqref{upperlower} it follows
\begin{align} \varphi(s) \le \frac{\psi(s)}{m} \le 2\varphi(s)
\label{upperlower2}
\end{align}
with
\[ \psi(s) := \frac{\sum_{y=1}^\infty f [y](y-1) (1+ s  + \cdots + s^{y-1})}{\sum_{z=1}^\infty f [z] (1+ s  + \cdots + s^{z-1})} \ .\]

Now consider the probability measures $g_s\in \mathcal P$, $0\le s \le 1$, given by
\begin{align} \label{gs} g_s[y] := \frac {f[y] (1+s+ \cdots + s^{y-1})}{ \sum_{z=1}^\infty f [z ] (1+ s  + \cdots + s^{z-1})} \ , \quad y \ge 1 \ .\end{align}
Then for $f[y] >0$ and $z>y$, after some algebra,
\[ \frac {g_s[z]}{g_s[y]} = \frac{f[z ]}{f[y ]} \prod_{v=1}^{z-y}\Big( 1+ \frac 1{s^{-1}+ \cdots + s^{-y-v+1}}\Big)\ , \]
which is an increasing function in $s$. Therefore by Lemma 3 the function $\psi(s)$ 
is increasing in $s$. In combination with \eqref{upperlower2} we get
\[ \varphi(s) \le \frac{\psi(s)}{m} \le \frac{\psi(1)}{m} \le 2 \varphi(1)\ , \ 
 2\varphi(s) \ge \frac{\psi(s)}{m}  \ge \frac{\psi(0)}{m} \ge  \varphi(0) \ . \]
 This gives the claim of the lemma.
\end{proof}

\begin{proof}[Proof of Lemma 2]
First we estimate the derivative of $\varphi$, which is given by
\[  \varphi'(s)= \frac 1m \frac{mf'(s)}{(1-f(s))^2} - \frac 1{m(1-s)^2} \ . \]
It turns out that this expression becomes more manageable if we replace the  squared  geometric mean $\sqrt{mf'(s)}$ on the right-hand side by the square of the arithmetic mean $(m+f'(s))/2$. Therefore we split the derivative into parts according to
\begin{align} \label{diff}  \varphi'(s) = \psi_1(s) - \psi_2(s)  \end{align}
with
\[ \psi_1(s)= \frac 1{4m} \frac{(m+f'(s))^2}{(1-f(s))^2} - \frac 1{m(1-s)^2} \ , \ \psi_2(s)= \frac 1{4 m}\frac{(m+f'(s))^2}{(1-f(s))^2} - \frac{f'(s)}{(1-f(s))^2} \ . \]
We show that both $\psi_1$ and $\psi_2$  are non-negative functions and estimate them from above.

For $\psi_1$ we accomplish this task by introducing the function
\begin{align*}
\zeta(s)&:= (m+f'(s))-2\frac{1-f(s)}{1-s}
\\
&= \sum_{y = 1}^\infty  y(1+s^{y-1})f[y]- 2\sum_{y = 1}^\infty\frac{1-s^y}{1-s} f[y]\\
& =\sum_{y = 3}^\infty \Big( y(1+s^{y-1}) - 2(1+ s+ \cdots + s^{y-1}) \Big) f[y] \ . 
\end{align*}
Since 
\begin{align*}\frac d{ds} \big(y(1+s^{y-1})& - 2(1+ s+ \cdots + s^{y-1})\big)\\&= y(y-1)s^{y-2} -2(1+2s + \ldots+ (y-1)s^{y-2})\\
&\le y(y-1)s^{y-2} -2s^{y-2}(1+2+ \ldots+ (y-1)) =0
\end{align*} 
for all $0 \le s \le 1$, and since $\zeta(1)=0$ we see that $\zeta$ is a non-negative, decreasing function. 
Thus $\psi_1$ is a non-negative function, too. 
Also $\zeta(0)\le m$. 

Moreover we have for $y \ge 3$ the polynomial identity  
\begin{align*}
y(1+s^{y-1}) - 2(1+ s+ \cdots + s^{y-1}) = (1-s)^2 \sum_{z=1}^{y-2} z(y-z-1)s^{z-1}  \ ,
\end{align*}
and consequently
\[\zeta(s)=  (1-s)^2 \xi(s) \]
with
\[ \xi(s):= \sum_{y=3}^\infty \sum_{z=1}^{y-2} z(y-z-1)s^{z-1}f[y] \ .\]
The function $\xi $ is  non-negative and increasing.  

Coming back to $\psi_1$ we rewrite it as
\[\psi_1(s)= \frac{\frac 12 (m+f'(s))(1-s) -(1-f(s))}{(1-f(s))(1-s)}\cdot\frac{\frac 12 (m+f'(s))(1-s) +(1-f(s))}{m(1-f(s))(1-s)}\ .\]
Using $f'(s)\le m$ it follows  
\begin{align*}
\psi_1(s)& \le \frac {\zeta(s)}{2(1-f(s))} \Big(\frac 1{1-f(s)} +\frac 1{m(1-s)})\Big)\\
& = \frac {\zeta(s)}{2 } \Big( \frac 1{m(1-s)} + \varphi(s)\Big)\Big( \frac 2{m(1-s)}+ \varphi(s)\Big)\\
& \le \ 2\zeta(s)  \Big( \frac 1{m^2(1-s)^2} + \varphi(s)^2 \Big) \\
&=  \frac{2 \xi(s)}{m^2} + 2 \zeta(s) \varphi(s)^2 \ .
\end{align*}
By means of  Lemma 1, by the monotonicity properties of $\xi$ and $\zeta$ and by $\varphi(1)=\nu/2$, $\zeta(0)\le m$ we obtain
\begin{align} 0\le  \psi_1(s) \le \frac{ 2\xi(s)}{m^2}  + 2m \nu^2 \ . \label{psi1} \end{align}

Now we investigate the function $\psi_2$, which we rewrite as
\[ \psi_2(s)= \frac 1{4m} \Big(\frac{m-f'(s)}{1-f(s)}\Big)^2 \ . \]
We have
\[ 1-f(s)= \sum_{z=1}^\infty (1-s^z)f[z]= (1-s)\sum_{z=1}^\infty(1+s+ \cdots + s^{z-1})f[z] \]
and
\[ m-f'(s)= \sum_{y=1}^\infty (1-s^{y-1})y f[y] = (1-s) \sum_{y=2}^\infty y(1+ \cdots + s^{y-2} )f[y] \ . \]
Using the notation from \eqref{gs} it follows
\[ \frac{m-f'(s)}{1-f(s)} = \sum_{y=2}^\infty \frac{1+ \cdots + s^{y-2}}{1+ \cdots + s^{y-1}} y g_s[y] \le\sum_{y=2}^\infty y g_s[y] \ .\]
As above we may apply Lemma 3 to the probability measures $g_s$ and conclude that the right-hand term is increasing with $s$. Therefore
\[ 0\le \frac{m-f'(s)}{1-f(s)} \le \sum_{y=2}^\infty y g_1[y] = \frac{\sum_{y=2}^\infty y^2 f[y]}{\sum_{z=1}^\infty z f[z]} \le \frac{2\sum_{y=1}^\infty y(y-1) f[y]}{\sum_{z=1}^\infty z f[z]}= 2m\nu \]
and hence
\begin{align} 0 \le \psi_2(s)\le m\nu^2 \ . \label{psi2}
\end{align}

Coming to our claim note first that owing to the non-negativity of $\psi_1$ and $\psi_2$ we obtain from formula \eqref{diff} for any $s \le u \le 1$
\begin{align*} - \int_s^1 \psi_2(t)\, dt \le \varphi(1)-\varphi(u) \le \int_s^1 \psi_1(t)\, dt \ . 
\end{align*}
The equations \eqref{psi1} and \eqref{psi2} entail
\begin{align} -m\nu^2 (1-s) \le \varphi(1)-\varphi(u) \le \frac 2{m^2}\int_s^1 \xi(t)\, dt + 2m \nu^2 (1-s)
\label{phiminusphi} \ .
\end{align}
It remains to estimate the right-hand integral. We have for $0 \le s < 1$
\begin{align*}\int_s^1\xi(t)\, dt &= \sum_{y=3}^\infty \sum_{z=1}^{y-2} (y-z-1)(1-s^{z})f[y] \\
& \le (1-s) \sum_{y=3}^\infty(y-2)f[y] \sum_{z=1}^{y-2} \sum_{u=0}^{z-1} s^u \\
&= (1-s) \sum_{y=3}^\infty(y-2)f[y] \sum_{u=0}^{y-3}(y-2-u)s^u \\
&= (1-s) \sum_{u=0}^\infty s^u \sum_{y=u+3}^\infty (y-2)^2 f[y] \ .
\end{align*}
The right-hand sum is monotonically decreasing in $u$, therefore for natural numbers $a$ we end up with the estimate
\begin{align*}
\int_s^1& \xi(t)\, dt \\& \le \sum_{y=3}^\infty (y-2)^2 f[y](1-s) \sum_{u=0}^{a-1} s^u  + \sum_{y=a+3}^\infty (y-2)^2 f[y](1-s) \sum_{u=a}^\infty s^u\\
&\le f''(1) a(1-s) + \mathbf E[Y^2; Y > a] \ .
\end{align*}
Combining this estimate with \eqref{phiminusphi} our claim follows.
\end{proof}

\paragraph{Remark 2.} We have
\[ \xi(1)= \sum_{y=3}^\infty \sum_{z=1}^{y-2} z(y-z-1)f[y]= \frac 13 \sum_{y=3}^\infty z(z-1)(z-2)f[z] = \frac {f'''(1)}3  \]
and hence from \eqref{diff}, \eqref{psi1}, \eqref{psi2} and the monotonicity of $\xi$ for $0 \le s \le 1$
\[ -\frac{f''(1)^2}{f'(1)^3} \le \varphi'(s) \le \frac{ 2f'''(1)}{3 f'(1)^2}  + 2 \frac{f''(1)^2}{f'(1)^3} . \]
The quality of these bounds becomes evident from the observation that 
\[ \varphi'(1) = \frac 16 \frac{f'''(1)}{f'(1)^2} - \frac 14 \frac {f''(1)^2}{f'(1)^3} \ , \]
as  follows by means of  Taylor expansions of $f$ and $f'$ about 1.
\qed

\section{Proof of the theorems}

First let us consider some formulas for moments. There exists a clear-cut expression for the variance of $Z_n$ due to Fearn \cite{Fea}. It seems to be less noticed that there is a similar appealing formula for the second factorial moment of $Z_n$, which turns out to be more useful for our purpose.

\paragraph{Lemma 4.} {\em For a BPVE $(Z_n)_{n \ge 0}$ we have}
\begin{align*} 
\mathbf E[Z_n] = \mu_n \ , \ \frac{\mathbf E[Z_n(Z_{n}-1)]}{
\mathbf E[Z_n]^2}=\sum_{k=1}^{n} \frac{\nu_{k}}{\mu_{k-1} } \ .
\end{align*}

\mbox{}\\
The proof follows a standard pattern. Let $v=(f_1,f_2, \ldots)$ denote a varying environment. 
For non-negative integers $k \le n$ let us define the probability measures 
\[ f_{k,n} := f_{k+1} \circ \cdots \circ f_n  \]
with the convention $f_{n,n}= \delta_1$ (the dirac measure at point 1). 
We have
\[ f_{k,n}'(s)= \prod_{l=k+1}^{n} f_l'(f_{l,n}(s)) \ , \]
in particular $f_{n,n}'(s)=1$, and after some rearrangements
\begin{align*}
 f_{k,n}''(s) = f_{k,n}'(s)^2\sum_{l=k+1}^{n} \frac{f_l''(f_{l,n}(s))}{f_l'(f_{l,n}(s))^2\prod_{j=k+1}^{l-1}f_j'(f_{j,n}(s))} \ ,  
 \end{align*}
in particular $f_{n,n}''(s)=0$. Since the distribution of $Z_n$ is given by $f_{0,n}$, choosing $k=0$ and $s=1$  Lemma 4 is proved.

Next we recall an expansion of the generating function of $Z_n$ taken from \cite{Ji} and \cite{Gei}. This kind of formula  has been used in many investigations of branching processes. Let  $\varphi_n$, $ n \ge 1$, be the shape functions of $f_n$, $n \ge 1$.
Then, since $f_{k,n}=f_{k+1}\circ f_{k+1,n}$ for $k <n$, we have
\[ \frac{ 1}{1-f_{k,n}(s)} = \frac{1}{f_{k+1}'(1)(1-f_{k+1,n}(s))} + \varphi_1(f_{k+1,n}(s)) \ .\]
Iterating the formula we end up with the following identity.

\paragraph{Lemma 5.} {\em For  $0\le s < 1$, $0 \le k < n$
\[ \frac{1}{1- f_{k,n}(s)} = \frac {\mu_k}{\mu_n(1-s)}+ \varphi_{k,n}(s) \quad \text{with} \quad \varphi_{k,n}(s):= \mu_k\sum_{l=k+1}^n \frac{\varphi_l(f_{l,n}(s))}{\mu_{l-1}}\ , \]
i.e. $\varphi_{k,n}$ is the shape function of $f_{k,n}$.}

\mbox{}\\
In order to estimate survival probabilities, assumtion \eqref{A1} now comes into play.
The next lemma reveals its role.

\paragraph{Lemma 6.} {\em Condition \eqref{A1} is fulfilled if and only if there is a constant $c'< \infty$ such that we have $\varphi_n(1)\le c'\varphi_n(0)$ for all $n \ge 1$.}

\begin{proof}
Recall that $Y_n$ denotes a random variable with distribution $f_n$. We have $\mathbf P(Y_n \ge 2)=0$ iff $\varphi_n(1)= \mathbf E[Y_n(Y_n-1)]/(2\mathbf E[Y_n]^2)=0$. Then both inequalities from \eqref{A1} and from our lemma are valid for all $c>0$ and $c'>0$, respectively. Therefore we may without loss of generality  assume that $\mathbf P(Y_n \ge 2)>0$ for all $n \ge 1$.
Then we have
\[ \varphi_n(0) = \frac 1{1-f_n[0]} - \frac 1{f_n'(1)} = \frac{ \mathbf E[(Y_n-1); Y_n \ge 1]}{\mathbf E[Y_n]\mathbf P(Y_n \ge1) } \]
and therefore because of \eqref{phi1}
\[  \frac{\varphi_n(1)}{\varphi_n(0)} = \frac{\mathbf E[Y_n(Y_n-1)]\mathbf P(Y_n \ge1) }{2\mathbf E[(Y_n-1); Y_n \ge 1]\mathbf E[Y_n] }\ .  \]
It is not difficult to see that these expressions are bounded uniformly in $n$ iff the same holds true for the terms
\[\frac{\mathbf E[Y_n^2; Y_n \ge 2  ]\mathbf P(Y_n \ge1)} { \mathbf E[Y_n; Y_n \ge 2] \mathbf E[Y_{n}  ] } \ ,\]
which in turn is equivalent to condition \eqref{A1}. This gives our claim.
\end{proof}

\noindent
In particular, if $\varphi_n(1)\le c'\varphi_n(0)$ for all $n \ge 1$ then we obtain for the shape functions $\varphi_{k,n}$ of the generating functions $f_{k,n}$ from Lemma 5 by means of Lemmas 6 an 1
\[  \varphi_{k,n}(1) = \mu_n\sum_{l=k+1}^n \frac{\varphi_l(1)}{\mu_{l-1}} \le c' \mu_n\sum_{l=k+1}^n \frac{\varphi_l(0)}{\mu_{l-1}}  \le 2c' \mu_n\sum_{l=k+1}^n \frac{\varphi_l(f_{l,n}(0))}{\mu_{l-1}} = 2c'\varphi_{k,n}(0)
\]
for all $1\le k \le n$. This estimate together with Lemma 6 prove our Remark 1 from the Introduction, namely that any subsequence of a regular BPVE is regular, too.

The next lemma has a forerunner in Agresti's estimate \cite[Theorem 1]{Agr}.

\paragraph{Lemma 7.} {\em Under Assumption \eqref{A1} there is a $\gamma>0$ such that for all $n \ge 0$}
\[  \frac{\mathbf E[Z_n]^2}{\mathbf E[Z_n^2]} \le \mathbf P(Z_n >0) \le \frac 1 \gamma \frac{\mathbf E[Z_n]^2}{\mathbf E[Z_n^2]}\ . \]

\begin{proof}
The left-hand estimate is just the standard Paley-Zygmund inequality. For the right-hand estimate observe that $\mathbf P(Z_n >0)= 1-  f_{0,n}[0]=1-f_{0,n}(0)$. Using Lemma 5 with $s=0$ we get the representation
\begin{align}
\frac 1{\mathbf P(Z_n >0)} = \frac 1{\mu_n} + \sum_{k=1}^n \frac {\varphi_k(f_{k,n}(0))}{\mu_{k-1} } \ , 
\label{representation}
\end{align}
hence by means of Lemma 1 
\begin{align}
\frac 1{\mathbf P(Z_n >0)}  \ge \frac 1{\mu_n} + \frac 12 \sum_{k=1}^n \frac {\varphi_k(0)}{\mu_{k-1} }\ . 
\label{Jir}
\end{align}
and by assumption \eqref{A1}, Lemma 6 and \eqref{phi1}
\[\frac 1{\mathbf P(Z_n >0)}  \ge \frac 1{\mu_n} + \frac 1{2c'} \sum_{k=1}^n \frac {\varphi_k(1)}{\mu_{k-1} } = \frac 1{\mu_n} + \frac 1{4c'} \sum_{k=1}^n \frac {\nu_k}{\mu_{k-1} } \ . \]
Letting $\gamma:= \min(1, (4c')^{-1})$ we obtain
\begin{align*} 
\frac 1{\mathbf P(Z_n >0)}  \ge \gamma \Big(\frac 1{\mu_n} + \sum_{k=1}^n \frac {\nu_k}{\mu_{k-1} }\Big) \ .
\end{align*}
On the other hand Lemma 4 implies
\begin{align} \label{secmom} \frac{\mathbf E[Z_n^2]}{\mathbf E[Z_n]^2} = \frac{\mathbf E[Z_n(Z_n-1)]}{\mathbf E[Z_n]^2} + \frac 1{\mathbf E[Z_n]} = \sum_{k=1}^n \frac{\nu_k}{\mu_{k-1}} + \frac 1{\mu_n} \ .\end{align}
Combining the last two formulas our claim follows.
\end{proof}

\begin{proof}[Proof of Theorem 1] 
(i) $\Leftrightarrow$ (ii):
Since $\lim_{n \to \infty} \mathbf P(Z_n>0) = 1-q$ the equivalence follows from  Lemma 7.

(ii) $\Leftrightarrow$ (iii): We have
\begin{align}
\sum_{k=1}^n \frac {\rho_k}{\mu_{k-1}}  &= \sum_{k=1}^n \frac{\nu_k+ f_k(1)^{-1}-1}{\mu_{k-1}}\notag \\&= \sum_{k=1}^n \frac {\nu_k}{\mu_{k-1}} + \sum_{k=1}^n \Big(\frac 1{\mu_k}-\frac 1{\mu_{k-1}}\Big) = \sum_{k=1}^n \frac {\nu_k}{\mu_{k-1}}  + \frac 1{\mu_n} -1 \ ,
\label{equiv}
\end{align}
thus because of \eqref{secmom}
\begin{align} \frac{\mathbf E[Z_n^2]}{\mathbf E[Z_n]^2} = \sum_{k=1}^n \frac {\rho_k}{\mu_{k-1}} +1 \ . 
\label{altern}
\end{align}
This gives the claim.

(iii) $\Leftrightarrow$ (iv): This equivalence is an immediate consequence of \eqref{equiv}.

(v) $\Leftrightarrow$ (vi): This implication follows again from  Lemma 7.

(vi) $\Leftrightarrow$ (vii): This is a consequence of equation \eqref{altern}.

(vii) $\Leftrightarrow$ (viii): Again this claim follows from \eqref{equiv}.
\end{proof}

\paragraph{Remark 3.} From \eqref{Jir} it follows that   a sufficient condition for a.s. extinction is given by the single requirement $\sum_{k\ge 1} \varphi_k(0)/\mu_{k-1} = \infty$ (without \eqref{A1}). This confirms a conjecture of Jirina~\cite{Ji}. \qed

\begin{proof}[Proof of Theorem 2]
Obviously statement (i) is valid.  For the first part of statement (ii) note that from Theorem 1, (vi) it follows that 
$\sup_{n \ge 0} \mathbf E[W_n^2]  < \infty$.
Therefore the martingale $(W_n)_{n \ge 0}$ is bounded in $\mathcal L^2$ implying $\mathbf E[W]=\mathbf E[W_0]=1$ and $\mathbf E[W^2]< \infty$. From \eqref{altern} it follows that
\[ \mathbf E[W^2] =  \sum_{k=1}^\infty \frac {\rho_k}{\mu_{k-1}} +1\ . \]
This implies formula \eqref{varW}.

For the proof of the last claim we distinguish two cases. Either $\mu_n \to r$ with $0<r<\infty$. Then $W_n =Z_n/\mu_n \to Z_\infty/r$ a.s., consequently $W=Z_\infty/r$ a.s. and $\mathbf P(W=0)= \mathbf P(Z_\infty=0)=q$.
Else we may assume $\mu_n \to \infty$  in view of Theorem 1, (viii). Also $ \{Z_\infty=0\} \subset \{W=0\}$ a.s., thus it is sufficient to show that $\mathbf P(Z_\infty >0, W=0)=0$.
First we estimate $\mathbf P(Z_\infty=0 \mid Z_k=1)$ from below.
From Lemma 5 and Lemma 1  for $k<n$
\[ \frac1{1- \mathbf P(Z_n=0\mid Z_k=1)} = \frac 1{1- f_{k,n}(0)} \ge \frac 12 \mu_k\sum_{l=k+1}^n \frac {\varphi_l(0)}{\mu_{l-1}} \ . \]
as well as  
\begin{align*}
\frac 1{1- \mathbf E[ e^{-\lambda W_n}\mid Z_k=1]} &= \frac 1{1-f_{k,n}(e^{-\lambda/\mu_n})} \\&\le \frac {\mu_k}{\mu_n(1- e^{-\lambda/\mu_n})} + 2 \mu_k\sum_{l=k+1}^n \frac {\varphi_l(1)}{\mu_{l-1}} 
\end{align*}
with $\lambda >0$. By means of Lemma 6 this entails
\[\frac 1{1- \mathbf E[ e^{-\lambda W_n}\mid Z_k=1]} \le \frac {\mu_k}{\mu_n(1- e^{-\lambda/\mu_n})} + \frac{4c'}{1- \mathbf P(Z_n=0\mid Z_k=1)} \ .\]
 Letting $n \to \infty$ we get
\[ \frac 1{1- \mathbf E[ e^{-\lambda W}\mid Z_k=1]} \le \frac {\mu_k}\lambda +\frac{4c'}{1- \mathbf P(Z_\infty=0\mid Z_k=1)} \]
and with $\lambda \to \infty$
\[ \frac 1{\mathbf P(W>0\mid Z_k=1)} \le \frac{4c'}{\mathbf P(Z_\infty>0 \mid Z_k=1)} \ . \]
Using $e^{-2x} \le 1-x $ for $0 \le x \le 1/2$ it follows for $\mathbf P(W>0 \mid Z_k=1) \le (8c')^{-1}$ that
\begin{align}
\mathbf P(Z_\infty&=0\mid Z_k=1)=1- \mathbf P(Z_\infty>0\mid Z_k=1)  \ge 1- 4c' \mathbf P(W>0\mid Z_k=1)\notag \\ &\ge e^{- 8c'\mathbf P(W>0\mid Z_k=1)} \ge (1- \mathbf P(W>0\mid Z_k=1))^{8c'}\notag \\& = \mathbf P(W=0\mid Z_k=1)^{8c'} \ . \label{Souza2}
\end{align}

Now we draw on a martingale, which already appears in the work of D'Souza and Biggins \cite{Bi}. Let for $n \ge 0$
\[ M_n := \mathbf P(W=0\mid Z_0, \ldots , Z_n) = \mathbf P(W=0\mid Z_n=1)^{Z_n} \text{ a.s.}  \]
From standard martingale theory $M_n\to I\{W=0\}$ a.s. In particular we have
\begin{align}\label{Souza} \mathbf P(W=0\mid Z_n=1)^{Z_n} \to 1  \text{ a.s. on the event that } W=0 \ ,
\end{align}
a result which has already been  exploited by D'Souza \cite{Sou}. 

We distinguish two cases. Either there is an infinite sequence of natural numbers such that $\mathbf P(W>0 \mid Z_n=1) > (8c')^{-1}$ along this sequence. Then \eqref{Souza} implies that $Z_n \to 0$ a.s. on the event $W=0$. Or else we may apply our estimate  \eqref{Souza2} to obtain from \eqref{Souza} that
\[ \mathbf P(Z_\infty=0 \mid Z_n=1)^{Z_n} \to 1 \text{ a.s. on the event that }W=0 \ . \]
Therefore, given $\varepsilon >0$, we have for $n $ sufficiently large
\begin{align*}
\mathbf P(Z_\infty >0, W=0) &\le \varepsilon + \mathbf P\big(Z_n >0, \mathbf P(Z_\infty=0 \mid Z_n=1)^{Z_n} \ge 1-\varepsilon\big) \\
& \le \varepsilon + \frac 1{1-\varepsilon} \mathbf E [ \mathbf P(Z_\infty=0\mid Z_n); Z_n>0] \\
&= \varepsilon + \frac 1{1-\varepsilon}  \mathbf P(Z_\infty=0, Z_n>0) \ . 
\end{align*}
Letting $n \to \infty$ we obtain $\mathbf P(Z_\infty >0, W=0) \le \varepsilon$, and the claim follows with $\varepsilon \to 0$.
\end{proof}

\begin{proof}[Proof of Theorem 3]
We begin with the proof of the last claim. Note that the assertion from Lemma 7 can be rewritten as 
\[ \gamma \frac {\mathbf E[Z_n^2]}{\mathbf E[Z_n]} \le \mathbf E[Z_n \mid Z_n>0] \le  \frac {\mathbf E[Z_n^2]}{\mathbf E[Z_n]} \]
and \eqref{secmom} gives
$  \mathbf E[Z_n^2]/\mathbf E[Z_n]= 1+ \mu_n \sum_{k=1}^n \frac{\nu_k}{\mu_{k-1}}=a_n$.
This implies \eqref{EZn}.

Consequently, by means of Markov's inequality we obtain
\[\mathbf P(Z_n/a_n >u \mid Z_n >0) \le \frac 1{ua_n} \mathbf E[Z_n \mid Z_n >0] \le \frac 1u \ ,\]
which implies the theorem's first claim.

Concerning the second claim we remark that for $a_n < 2$ we may set $u=1/2$. For $a_n \ge 2$ we have by means of Lemma 5 the estimate
\begin{align*} 1-s^u + \mathbf P(Z_n/a_n >u) &\ge \mathbf E[1- s^{Z_n/a_n}\mid Z_n >0] \mathbf E[ 1-s^{Z_n/a_n}\mid Z_n>0]= \frac{1-f_{0,n}(s^{1/a_n})}{1-f_{0,n}(0)}
\\&= \Big(\frac 1{\mu_n} + \sum_{k=1}^n \frac{\varphi_k(f_{k,n}(0))}{\mu_{k-1}}\Big)\Big/\Big(\frac 1{\mu_n(1-s^{1/a_n})} + \sum_{k=1}^n \frac{\varphi_k(f_{k,n}(s^{1/a_n}))}{\mu_{k-1}}\Big) 
\end{align*}
with $0<s<1$ and $u>0$.  Lemma 1, Lemma 6 and \eqref{phi1} yield the bound
\begin{align*} 1-s^u + \mathbf P(Z_n/a_n >u) &\ge  \sum_{k=1}^n \frac{\varphi_k(0)}{2\mu_{k-1}}\Big/\Big(\frac 1{\mu_n(1-s^{1/a_n})} + 2\sum_{k=1}^n \frac{\varphi_k(1)}{\mu_{k-1}}\Big)\\ &\ge
\frac 1{4c'}\sum_{k=1}^n \frac{\nu_k}{\mu_{k-1}}\Big/\Big(\frac 1{\mu_n(1-s^{1/a_n})} + \sum_{k=1}^n \frac{\nu_k}{\mu_{k-1}}\Big)
\end{align*}
Moreover $1-s^{1/a_n} \ge  a_n^{-1}(1-s) $, since $1/a_n \le 1$. Hence, choosing $s=1/2$ we get
\[1-2^{-u} + \mathbf P(Z_n/a_n >u) \ge \frac 1{4c'}\sum_{k=1}^n \frac{\nu_k}{\mu_{k-1}}\Big/\Big(\frac {2a_n}{\mu_n} + \sum_{k=1}^n \frac{\nu_k}{\mu_{k-1}}\Big) \]
Finally, from $a_n\ge 2$ it follows that $a_n \le  2 \mu_n \sum_{k=1}^{n} \nu_k/\mu_{k-1}$ and consequently
\[ 1-2^{-u} + \mathbf P(Z_n/a_n >u) \ge \frac 1{20c'}  \]
for all $u>0$. If we set now $\theta =1/(40 c')$ and choose $u>0$ so small that $1-2^{-u} \le \theta$ we obtain
$ \mathbf P(Z_n/a_n >u) \ge \theta $
which is our second claim.
\end{proof}

\mbox{}\\
The next lemma prepares the proof of Theorem 4. It clarifies the role of condition \eqref{B}.

\paragraph{Lemma 8.} {\em Assume condition \eqref{B} and  let $q=1$. Then the condition
$ 1/{\mu_n}= o\big(\sum_{k=1}^n  \nu_k/\mu_{k-1}\big)$
implies
\begin{align*} \sup_{0 \le s \le 1}\Big| \sum_{k=1}^n\frac{ \varphi_k(f_{k,n}(s))}{\mu_{k-1}} -   \sum_{k=1}^n\frac{\varphi_k(1)}{\mu_{k-1}} \Big| = o\Big(\sum_{k=1}^n\frac{ \varphi_k(1)}{\mu_{k-1}}\Big)
\end{align*}
as $n\to \infty$.}

\begin{proof} 
Fix $\varepsilon >0$ and choose $c_{\varepsilon/9}$ according to assumption \eqref{B}. Let
\[ s_k := 1- \frac \eta{1+ f_k'(1)} \]
with some $0<\eta < 1$. Then from Lemma 2 with $a= \lfloor c_{\varepsilon/9} \rfloor$
\begin{align*}
\sup_{s_k \le t \le 1}|\varphi_k(1)-\varphi_k(t)| \le  2\nu_k\frac{ f_k''(1) }{f_k'(1)} \frac \eta{ 1+ f_k'(1)} +2 c_{\varepsilon/9 } \nu_k \eta + \frac \varepsilon 9  4\nu_k \ .
\end{align*}
From the estimate \eqref{BimpliesA} it follows that 
\begin{align}
f_k''(1) \le 2c_{1/2}f_k'(1)(1+ f_k'(1))\ .
\label{B2}
\end{align}
Therefore there is a $\eta=\eta_\varepsilon>0$ such that
\begin{align}  \sup_{s_k \le t \le 1}|\varphi_k(1)-\varphi_k(t)| \le \frac \varepsilon 2 \nu_k =\varepsilon \varphi_k(1) \ . 
\label{sup}
\end{align}

Now set
\[ r=r_{\varepsilon,n} :=  \min\{k \le n : f_{k,n}(0)  \le s_k\} \ . \]
Because of $f_{n,n}(0)=0$ this minimum is attained. In view of \eqref{sup} and Lemma 1 it follows
\begin{align*}
 \Big| \sum_{k=1}^n\frac{ \varphi_k(1)}{\mu_{k-1}} -\sum_{k=1}^n\frac{ \varphi_k(f_{k,n}(s))}{\mu_{k-1}} \Big| \le  \varepsilon  \sum_{k=1}^{r-1}\frac{ \varphi_k(1)}{\mu_{k-1}}  +3 \frac{\varphi_r(1)}{\mu_{r-1}}+ 3\sum_{k=r+1}^n\frac{ \varphi_k(1)}{\mu_{k-1}} \ .
\end{align*}
From \eqref{B2} we have
\begin{align*}
\frac{\varphi_r(1)}{\mu_{r-1}}= \frac{f_r''(1)}{2f_r'(1)^2 \mu_{r-1}} \le \frac{c_{1/2}( f_r'(1)+1)}{f_r'(1)\mu_{r-1}} = c_{1/2} \Big(  \frac 1{\mu_{r-1}}+\frac 1{\mu_r}\Big)  
\end{align*}
and from Lemma 6
\begin{align*}
\sum_{k=r+1}^n \frac{ \varphi_k(1)}{\mu_{k-1}} &\le c' \sum_{k=r+1}^n \frac{ \varphi_k(0)}{\mu_{k-1}} \le
2 c'  \sum_{k=r+1}^n \frac{ \varphi_k(f_{k,n}(0))}{\mu_{k-1}} \ .
\end{align*} 
From \eqref{representation} it follows that $\mathbf P(Z_n>0 \mid Z_r=1)^{-1} =\mu_r/\mu_n+ \mu_r \sum_{k=r+1}^n\varphi_k(f_{k,n}(0))/\mu_{k-1}$ for $n>r$, hence we may proceed to
\begin{align*}
\sum_{k=r+1}^n \frac{ \varphi_k(1)}{\mu_{k-1}} \le \frac{2c'}{\mu_r(1-f_{r,n}(0))} \le \frac {2c'}{\mu_r (1-s_r)} = \frac{2c'(f_r'(1)+1)}{\eta\mu_r} = \frac{2c'}{\eta}\Big( \frac 1{\mu_{r-1}}+\frac 1{\mu_r}\Big) \ .
\end{align*} 
Putting our estimates together we get
\begin{align}
 \Big| \sum_{k=1}^n\frac{ \varphi_k(1)}{\mu_{k-1}} -\sum_{k=1}^n\frac{ \varphi_k(f_{k,n}(s))}{\mu_{k-1}} \Big| \le  \varepsilon  \sum_{k=1}^{n}\frac{ \varphi_k(1)}{\mu_{k-1}} + 3\Big( c_{1/2} +\frac{2c'}{\eta}\Big) \Big(  \frac 1{\mu_{r-1}}+\frac 1{\mu_r}\Big) \ .
 \label{assumption2}
\end{align}

Now the assumption $1/\mu_n=o\big(\sum_{k=1}^{n} \nu_k/\mu_{k-1}\big)$ comes into play. It implies that there is a positive integer $r_\varepsilon$ such that for all $r,n$ with $ r_\varepsilon<r\le n $
\begin{align}3\Big( \frac{2c'}{\eta}+c_{1/2} \Big) \Big(  \frac 1{\mu_{r-1}}+\frac 1{\mu_r}\Big) \le \frac \varepsilon 2  \sum_{k=1}^{r-1}\frac{ \nu_k}{\mu_{k-1}}+\frac \varepsilon 2  \sum_{k=1}^{r}\frac{ \nu_k}{\mu_{k-1}} \le  \varepsilon  \sum_{k=1}^{n}\frac{ \varphi_k(1)}{\mu_{k-1}} \ . 
\label{assumption}
\end{align}
Also from the   assumptions $q=1$ and  $1/\mu_n=o\big(\sum_{k=1}^{n} \nu_k/\mu_{k-1}\big)$ together with  Theorem 1 (iv) and \eqref{phi1} we have
\[ \sum_{k=1}^n\frac{ \varphi_k(1)}{\mu_{k-1}}=\frac 12 \sum_{k=1}^n\frac{ \nu_k}{\mu_{k-1}} \to \infty  \]
as $n \to \infty$, which implies that \eqref{assumption} hold for all $r \le r_\varepsilon$ and thus for all $r \le n$, if only $n$ is large enough. Thereby we may combine \eqref{assumption2} and \eqref{assumption} and obtain
\[\Big| \sum_{k=1}^n\frac{ \varphi_k(f_{k,n}(s))}{\mu_{k-1}} -  \sum_{k=1}^n\frac{ \varphi_k(1)}{\mu_{k-1}} \Big| \le 2 \varepsilon  \sum_{k=1}^{n}\frac{ \varphi_k(1)}{\mu_{k-1}} \]
for sufficiently large $n$.
This proves our claim.
\end{proof}

\begin{proof}[Proof of Theorem 4] 

(i) $\Rightarrow$ (ii): We argue by contradiction. If  assertion  (ii) fails, then there is an increasing sequence $(n_i)_{i\ge 0}$ of natural numbers fulfilling $\sup_i \mathbf E[Z_{n_i}\mid Z_{n_i}>0] < \infty$. From Theorem 3 it follows that  the random variables $Z_{n_i}$, $i \ge 0$, conditioned on $Z_{n_i}>0$ are tight. This does not conform with assertion (i), which proves the implication.

(ii) $\Rightarrow$ (iii): This implication follows from Theorem 3, since the assertion from (iii)  just states  that $a_n\to \infty$.

(iii) $\Rightarrow$ (i):
For the proof let
\[ b_n:= \frac{\mu_n}2 \sum_{k=1}^n \frac{\nu_k}{\mu_{k-1}} \ .\]
From Lemma 5 we have
\begin{align*} 1- \mathbf E[ &e^{-\lambda Z_n/b_n}\mid Z_n>0] = \frac{1-f_{0,n}(e^{-\lambda/b_n})}{1-f_{0,n}(0)} \\&= 
\Big(\frac 1{\mu_n} + \sum_{k=1}^n \frac{\varphi_k(f_{k,n}(0))}{\mu_{k-1}}\Big)
\Big(\frac 1{\mu_n(1- e^{-\lambda/b_n})} + 
\sum_{k=1}^n \frac{\varphi_k(f_{k,n}(e^{-\lambda/b_n}))}{\mu_{k-1}}
\Big)^{-1} \ .
\end{align*}
Since $b_n \to \infty$, from Lemma 8 and the theorem's assumption
 \[ 1- \mathbf E[ e^{-\lambda Z_n/b_n}\mid Z_n >0]=
\Big((1+o(1))\sum_{k=1}^n \frac{\nu_k}{2\mu_{k-1}}\Big)
\Big((1+o(1))\frac {b_n}{\lambda\mu_n} + 
(1+o(1))\sum_{k=1}^n \frac{\nu_k}{2\mu_{k-1}}
\Big)^{-1}
\]
as $n \to \infty$. 

\newpage

From the definition of $b_n$ we get
\[1- \mathbf E[ e^{-\lambda Z_n/b_n}\mid Z_n>0] = \frac{\lambda+o(1)}{1+\lambda } \ .\]
This implies assertion (i).

Moreover, from \eqref{representation}, Lemma 8 and assertion (iii) it follows that
\[\frac 1{\mathbf P(Z_n >0)} = \frac 1{\mu_n} + \sum_{k=1}^n \frac {\varphi_k(f_{k,n}(0))}{\mu_{k-1} } \sim \frac 12 \sum_{k=1}^n \frac{\nu_k}{\mu_{k-1}} \ . \]
This formula give the extra claims, which concludes the proof.
\end{proof}

\begin{proof}[Proof of Proposition 1]
By Theorem 1 (viii) the condition $q<1$ is equivalent to the requirements of both $\sum_{k=1}^\infty \nu_k/\mu_{k-1}< \infty$ and $0< \lim_n \mu_n \le \infty$. As already explained, the division between the supercritical regime and the asymptotically non-degenerate regime corresponds to the cases $\lim_n \mu_n=\infty $ and $0< \lim_n \mu_n<\infty$. This gives the first two assertions of the proposition.

Next the critical regime is given by the requirements that both $\mathbf  E[Z_n \mid Z_n>0] \to \infty$ and $q=1$. By Theorem 3 and Theorem 1 (iv) we may equivalently require that $1/\mu_n =o(\sum_{k=1}^n \nu_k/\mu_{k-1})$ together with either $\sum_{k=1}^n \nu_k/\mu_{k-1}=\infty$ or $\mu_n \to 0$. However, the third and the first of these conditions imply the second one, there for the third condition can be skipped, and we end up with the requirements $1/\mu_n =o(\sum_{k=1}^n \nu_k/\mu_{k-1})$ and $\sum_{k=1}^n \nu_k/\mu_{k-1}=\infty$, as stated in the proposition.

Finally, the subcritical regime is characterized by the conditions $\mathbf  E[Z_n \mid Z_n>0] \not\to \infty$ and $q=1$. Because of Theorem 3 the first condition is equivalent to the requirement \mbox{$a_n \not\to \infty$} respectively to $\liminf_n \mu_n \sum_{k=1}^n \nu_k/\mu_{k-1} < \infty$. Moreover, $\liminf_n \mu_n=0$ implies $q=1$, therefore the conditions stated in the proposition imply subcriticality. Conversely, if $q=1$ then by \mbox{Theorem 1 (iv)} we have $\lim_n \mu_n=0$ or $\sum_{k=1}^n \nu_k/\mu_{k-1} < \infty$. The former of these conditions trivially yields $\liminf_n \mu_n=0$, whereas the latter 
together with $\liminf_n \mu_n \sum_{k=1}^n \nu_k/\mu_{k-1} < \infty$ implies  $\liminf_n \mu_n=0$. Therefore the two conditions stated in the proposition are as well necessary for subcriticality.
\end{proof}


\begin{thebibliography}{99}

\bibitem{Agr} A. Agresti, On the extinction times of random and varying environment
branching processes. {\em J. Appl. Probab.} {\bf 12} (1975), 39--46.

\bibitem{Ba} V. Bansaye, F. Simatos, On the scaling limit of Galton Watson processes in varying environment. {\em Electron. J. Probab.} {\bf 20} (2015), 36 pp.

\bibitem{Bha} N. Bhattacharya, M. Perlman, Time-inhomogeneous branching processes conditioned on non-extinction. Preprint (2017). arXiv:1703.00337 [math.PR]

\bibitem{BrHaut} P. Braunsteins, S. Hautphenne, Extinction in lower Hessenberg branching processes with countably many types. {\em Ann. Appl. Probab.} {\bf 29} (2019), 2782--2818.

\bibitem{Chu} J. D. Church, On infinite composition products of probability generating functions. {\em Z. Wahrscheinlichkeitstheorie verw. Geb.} {\bf 19} (1971), 243--256.

\bibitem{Sou} J. C. D'Souza, The rates of growth of the Galton-Watson process in varying environments. {\em Adv. Appl. Probab.} Ê{\bf 26} (1994), 698--714.

\bibitem{Bi}  J. C. D'Souza, J. D. Biggins, The supercritical Galton-Watson process in varying environments. {\em Stoch. Proc. Appl.} {\bf 42} (1992), 39--47.

\bibitem{Do} D. Dolgopyat, P. Hebbar, L. Koralov, M. Perlman, Multi-type branching processes with time-dependent branching rates.  {\em J. Appl. Probab.} {\bf 55} (2018), 701--727. 

\bibitem{Fah} K. S. Fahady, M. P. Quine, D. Vere Jones,  Heavy traffic approximations for the Galton-Watson process. {\em Adv. Appl. Probab}. {\bf 3} (1971), 282--300.

\bibitem{Fea}  D. H. Fearn, Galton-Watson processes with generation dependence. {\em Proc. 6th Berkeley Symp. Math. Statist. Probab.} {\bf 4} (1971), 159--172.

\bibitem{Gei} J. Geiger, G. Kersting, The survival probability of a critical branching process in random environment. {\em Theor. Probab. Appl.} {\bf  45}  (2001), 517--525.

\bibitem{Goe} R. T. Goettge, Limit theorems for the supercritical Galton-Watson process in varying environments.
{\em Math. Biosci.} {\bf  28} (1976), 171--190.

\bibitem{GoKeMiPu}  M. Gonz\'alez, G. Kersting, C. Minuesa, I. del Puerto, Branching processes in varying environment with generation dependent immigration. {\em Stochastic Models} {\bf 35} (2019), 148--166.

\bibitem{KeVa} G. Kersting, V.  Vatutin, {\em Discrete time Branching Processes in Random Environment}. Wiley, 2017.

\bibitem{Ja} P. Jagers, Galton-Watson processes in varying environments. {\em J. Appl. Probab.} {\bf 11} (1974), 174--178.

\bibitem{Ji} M. Jirina, Extinction of non-homogeneous Galton-Watson processes. {\em J. Appl. Probab.} {\bf 13} (1976),  132--137.

\bibitem{Li} T. Lindvall, Almost sure convergence of branching processes in varying and random environments. {\em Ann. Probab.} {\bf 2} (1974), 344--346.


\bibitem{Ru} R. Lyons, Random walks, capacity and percolation on trees. {\em Ann. Probab.} {\bf 20}  (1992), 2043--2088.

\bibitem{Schuh} I. M. MacPhee,    H. J. Schuh,  A Galton-Watson branching process in varying environments with essentially constant means and two rates of growth. {\em Austral. J. Statist.} {\bf 25} (1983), 329--338.


\bibitem{SaJa} S. Sagitov, J. Jagers, Rank-dependent Galton-Watson processes and their pathwise duals. {\em J. Appl. Probab.} {\bf 50(A)} (2019), 229--239.

\bibitem{Sev} B. A. Sevast'yanov,   Transient phenomena in branching stochastic processes. {\em Theor. Probab. Appl.} {\bf 4} (1959), 113--128.

\end{thebibliography}
\end{document}